\documentclass[12pt,a4paper]{amsart}
\usepackage{amsmath}
\usepackage{amscd}
\numberwithin{equation}{section}

\input{prepictex}
\input{pictex}
\input{postpictex}
\chardef\bslash=`\\ 

\makeatletter
\def\verbatim{\interlinepenalty\@M \@verbatim
  \leftskip\@totalleftmargin\advance\leftskip2pc
  \frenchspacing\@vobeyspaces \@xverbatim}
\makeatother
\newtheorem{theorem}{Theorem}[section]
\newtheorem{corollary}[theorem]{Corollary}
\newtheorem{lemma}[theorem]{Lemma}

\theoremstyle{definition}
\newtheorem{definition}[theorem]{Definition}
\newtheorem{remark}[theorem]{Remark}

\newtheorem{example}[theorem]{Example}
\theoremstyle{main}
\newtheorem{main}{Theorem}

\newcounter{picture}

\newsymbol \rtimes 226F 
\DeclareMathOperator{\conv}{conv}
\DeclareMathOperator{\Dim}{dim}

\DeclareMathOperator{\Star}{star}

\newcommand{\CC}{{\mathbb C}}
\newcommand{\FF}{{\mathbb F}}

\newcommand{\PP}{{\mathbb P}}
\newcommand{\QQ}{{\mathbb Q}}
\newcommand{\RR}{{\mathbb R}}

\newcommand{\ZZ}{{\mathbb Z}}

\newcommand{\cO}{{\mathcal O}}

\newcommand{\D}{{\Delta}}
\newcommand{\bD}{{\partial \D}}
\newcommand{\dz}{{{\delta}^0}}
\newcommand{\di}{{{\delta}^1}}
\newcommand{\dii}{{{\delta}^n}}
\newcommand{\dd}{{{\delta}^2}}
\newcommand{\e}{{\varepsilon}}
\newcommand{\G}{{\Gamma}}
\newcommand{\Om}{{\Omega}}
\newcommand{\g}{{\gamma}}
\newcommand{\s}{{\sigma}}
\newcommand{\w}{{\omega}}

\newcommand{\0}{{\bf 0}}
\newcommand{\1}{{\bf 1}}
\newcommand{\2}{{\bf 2}}
\newcommand{\p}{{\bf p}}
\newcommand{\ka}{{\bf k}}
\newcommand{\n}{{\bf n}}

\newcommand{\fD}{{\mathfrak D}}

\newcommand{\fS}{{\mathfrak S}}

\newcommand{\od}{{\overline d}}

\newcommand{\ove}{{\overline e}}
\newcommand{\he}{{\hat e}}
\newcommand{\op}{{\overline p}}
\newcommand{\oq}{{\overline q}}
\newcommand{\hq}{{\hat q}}

\newcommand{\tA}{{\widetilde A}}

\newcommand{\zomatrix}{{$\{0,1\}$-matrix}}

\newcommand{\id}{{\bf I}}
\newcommand{\aut}{{\text{\rm Aut}}}

\newcommand{\PGL}{{\text{\rm{PGL}}}}
\newcommand{\PSL}{{\text{\rm{PSL}}}}

\newcommand{\wt}{\widetilde}

\begin{document}

\title[]{Torsion in $K$-theory for boundary actions on affine buildings of type $\tA_n$}

\date{November 17, 2000}
\author{Guyan Robertson}
\address{Mathematics Department, University of Newcastle, Callaghan, NSW 
2308, Australia}
\email{guyan@maths.newcastle.edu.au}
\thanks{This research was supported by the Australian Research Council and the Universit\'e d'Orl\'eans} 
\thanks{ \hfill Typeset by  \AmS-\LaTeX}

\begin{abstract}
Let $\G$ be a torsion free lattice in $G = \PGL(n+1,\FF)$, where $n\ge 1$ and $\FF$ is a non-archimedean local field.
Then $\G$ acts on the Furstenberg boundary $G/P$, where $P$ is a minimal parabolic subgroup of $G$. The identity element $\id$ in the crossed product $C^*$-algebra $C(G/P)\rtimes \G$ generates a class $[\id]$ in the $K_0$ group of $C(G/P)\rtimes \G$. It is shown that $[\id]$ is a torsion element of $K_0$ and there is an explicit bound for the order of $[\id]$. The result is proved more generally for groups acting on affine buildings of type $\tA_n$.
For $n=1, 2$ the Euler-Poincar\'e characteristic $\chi(\G)$ annihilates the class $[\id]$.   
\end{abstract}      

\maketitle

\section*{Introduction}

A. Connes \cite[Corollary 6.7]{connes}  proved that if $\G$ is a torsion free cocompact lattice in $\PSL (2,\RR)$ then the class of the identity
$[\id]$ in $K_0(C(\PP_1(\RR))\rtimes \G)$ has torsion. See also \cite{natsume,connes2,hn,ad}. More precisely, the order of $[\id]$ is  $-\chi(\G)$, where $\chi(\G)$ is the Euler-Poincar\'e characteristic of the corresponding Riemann surface \cite[Proposition 2.9]{ad}, \cite[Section 10]{hn}.
If $\G$ is a torsion free non-cocompact lattice in $\PSL (2,\RR)$, which means that $\G$ is a free group, then $[\id]=0$ in $K_0(C(\PP_1(\RR))\rtimes \G)$ \cite[Theorem 2.1]{ad}.
Moreover, Connes proved \cite[Corollary 5.9]{connes} that if
$\G$ is any countable subgroup of $\PSL (2, \CC)$ then the canonical map $K_0(C(\PP_1(\CC))) \to K_0(C(\PP_1(\CC))\rtimes \G)$ is injective and so in this case $[\id]$ is a non-torsion element in $K_0$. There is therefore a dramatic difference between the Fuchsian and the Kleinian cases.

Now let $\FF$ be a non-archimedean local field and let $\G$ be a torsion free (cocompact) lattice in $G = \PGL(2,\FF)$. As in the Fuchsian case, $\chi(\G).[\id]=0$ in $K_0(C(\PP_1(\FF))\rtimes \G)$. This assertion has a simple geometric proof, which is presented in Section \ref{section1} below. The Bruhat-Tits building associated with $\PGL(2,\FF)$ is a regular tree $\D$ whose boundary $\bD$ may be identified with the projective line $\PP_1(\FF)$. The result is proved in the more general situation where $\G$ is a discrete group  acting freely with compact quotient on  a tree $\D$. Such a group $\G$ is necessarily a free group \cite [Chapitre I.3.3]{ser}. 

A tree is an affine building of type $\tA_1$. The principal object of this paper is to prove corresponding results for group actions on buildings of type $\tA_n$, where $n\ge 2$. These buildings are in some respects higher dimensional analogues of trees, but have a very much more rigid structure.

 Let $\D$ be a building of type $\tA_n$. Then the link of any fixed vertex in $\D$ is the flag complex of a finite projective geometry. The order $q$ of this projective geometry is independent of the vertex, and is referred to as the order of the building $\D$. 
The case $n=2$ is particularly interesting because then the link of a vertex is a projective plane, which is possibly non-Desarguesian, and
$\D$ need not be the Bruhat-Tits building of a linear group.
The boundary $\bD$ is the set of chambers of the spherical building at infinity \cite[Chapter 9]{ron}, endowed with a natural totally disconnected compact Hausdorff topology \cite[Section 4]{ca2}.
Denote by $\aut(\D)$ the automorphism group of $\D$, equipped with the compact open topology. Then $\aut(\D)$ acts in a natural way on $\bD$.

Given an~$\tA_n$ building~$\D$ with vertex set $\D^0$, there is a type map
$\tau : \D^0 \to \ZZ/(n+1)\ZZ$ such that each maximal simplex ({\em chamber}) has exactly one vertex of each type.
An automorphism~$\alpha$ of~$\D$ is said to be {\em type-rotating} if
there exists $i\in \ZZ/(n+1)\ZZ$ such that $\tau(\alpha v)=\tau(v)+i$ for all
vertices $v\in\D^0$.

 The main result is the following.

\begin{main}\label{main}
Let $\G$ be a torsion free discrete group  of type rotating automorphisms of
a locally finite affine building $\D$ of type $\tA_n$. 
Suppose that $\G$ acts cocompactly on $\D$ with a finite number $n_0$ of vertex orbits. 
In the group $K_0(C(\bD)\rtimes \G)$, let $[\id]$ denote the class of the identity element.
\begin{itemize}
\item[(1)] If $n=1 \text{\ or\ } 2$, then $\chi(\G).[\id]=0$, where $\chi(\G)$ is the Euler-Poincar\'e characteristic of the group $\G$. \label{main1}
\item[(2)]  \label{main2}
If $n\ge 2$, then the order of $\D$ is a prime power $q$ and
$m.[\id]=0$, where
\begin{equation*}
m = 
\begin{cases}
    n_0(q-1)&  \text{if $n$ is odd},\\
    n_0(q^2-1) &  \text{if $n$ is even.}
\end{cases}
\end{equation*}
\end{itemize}
\end{main}
If $n=2$ then both parts of Theorem \ref{main} apply and can be combined to sharpen the result. This is done in Section \ref{EPC}.

The condition in Theorem \ref{main} that the building $\D$ is locally finite is an easy consequence of the other hypotheses, by Lemma \ref {inversionfree} below.
Also the finite cell complex $\G\backslash\D$ is a $K(\G,1)$ space and the group $\G$ has finite cohomological dimension \cite[I.1.5]{ser2}. It follows that  $\chi(\G)$ coincides with the usual Euler-Poincar\'e characteristic of the cell complex $\G\backslash\D$ \cite[I.1.5, Proposition 9]{ser2}.
Moreover $n_0$ is the number of vertices of $\G\backslash\D$. 
 
Let $\FF$ be a non-archimedean local field with residue field of order $q$.
The Bruhat-Tits building $\D$ of $G=\PGL(n+1,\FF)$ is a building of type $\tA_n$ and the action of $G$ on $\D$ is type rotating \cite{st}. 
The vertex set of $\D$ may be identified with the homogeneous space $G/K$, where $K$ is an open maximal compact subgroup of $G$. Then $n_0$ is the number of double cosets in $\G\backslash G/K$.
The boundary $\bD$ may be identified with the Furstenberg boundary $G/P$, where $P$ is a minimal parabolic subgroup of $G$. A torsion free lattice $\G$ in $G=\PGL(n+1,\FF)$ is automatically cocompact \cite[Chapitre II.1.5, p. 116]{ser}. The result therefore applies to such lattices, acting on the Furstenberg boundary and shows that the class $[\id]$ in $K_0(C(G/P)\rtimes \G)$ has torsion. 
If the Haar measure $\mu$ on $G$ is normalized so that $\mu(K)=1$, then $n_0$ is the covolume of $\G$ in $G$.
It follows from Theorem \ref{main} that there is a bound for the order of $[\id]$, which depends only on the covolume of $\G$ in $G$.

  The tree result is proved in Section~\ref{section1}, followed by the higher dimensional cases in Sections~\ref{section2} and~\ref{section3}. Section~\ref{EPC} is concerned with lattices in $\PGL(n+1,\FF)$ and computation of the Euler-Poincar\'e characteristic.

It is natural to ask if there is a similar result for lattices in the group $\PGL(n+1,\RR)$. The methods used here do not throw any light on this question.

It is worth remarking that the crossed product $C^*$-algebras considered here are
purely infinite, simple, separable, unital, nuclear and satisfy the Universal Coefficient Theorem. This was proved in \cite{rs2} for buildings of type $\tA_2$, and follows from \cite{jr} in higher dimensions. These algebras are therefore classified up to isomorphism by their two K-groups together with the class of $[\id]$ in $K_0$ \cite{k}. Properties of the class $[\id]$ are therefore particularly significant.

The author is grateful to the referees of this article for their helpful comments.

\section{The $\tA_1$ result}\label{section1}

The tree case is treated first, since the proof is relatively simple and motivates the harder $\tA_2$ case.  Two vertices of a tree have the same type if and only if the graph distance between them is even. Any automorphism of a tree is type rotating, so the type rotating hypothesis in Theorem~\ref{main} is vacuous in the $\tA_1$ case. 

Let $\D$ be a locally finite tree whose vertices all have degree at least three. The terminology of \cite{ser} will be used where possible. 
In particular the edges of $\D$ are directed edges. Each geometric edge of $\D$ corresponds to two directed edges $d$ and $\od$.
Let $\D^0$ denote the set of vertices and $\D^1$  the set of directed edges of $\D$.

Let $\G$ be a torsion free discrete subgroup of $\aut(\D)$ which acts cocompactly on $\D$. 
Then $\G$ acts freely (c.f. Lemma \ref{inversionfree} below), and so is a free group \cite[I.3.3]{ser}. In particular $\G$ acts without inversion, which means that no element $\g\in \G\setminus\{1\}$ stabilizes a geometric edge of $\D$. It follows that there is an orientation on the edges which is invariant under $\G$ \cite[I.3.1]{ser}. Choose such an orientation.  
This orientation consists of a partition $\D^1=\D^1_+\sqcup \overline{\D^1_+}$ and a bijective involution $d\mapsto \overline d : \D^1 \to \D^1$ which interchanges the two components of $\D^1$. Each directed edge $d$ has an
origin $o(d)\in \D^0$ and a terminal vertex $t(d)\in \D^0$ such that $o(\overline d)=t(d)$. The quotient $X=\G\backslash\D$ is a finite connected graph with vertex set $V=\G\backslash\D^0$
and directed edge set $E=E_+\sqcup\overline{E_+}$, where $E_+=\G\backslash\D^1_+$ and $\overline{E_+}=\G\backslash\overline{\D^1_+}$.  The graph $X$ has an induced involution $e\mapsto \overline e : E \to E$ and there are maps $o,t : E \to V$ with $o(\overline e)=t(e)$.
The Euler-Poincar\'e characteristic of the graph is $\chi(X)=n_0-n_1$ where $n_0=\#(V)$ and $n_1=\#(E_+)$.  

The boundary $\bD$ of the tree $\D$ can be identified with the set of equivalence classes of infinite semi-geodesics in $\D$, where equivalent semi-geodesics contain a common sub-semi-geodesic. Also $\bD$ has a natural compact totally disconnected topology \cite[I.2.2]{ser}. 

The group $\G$ acts on $\bD$ and hence on $C(\bD)$ via $\g\mapsto \alpha_{\g}$, where 
$\alpha_{\g} f(\w)=f(\g^{-1}\w)$, for $f\in C(\bD)$, $\g\in\G$.
The algebraic crossed product relative to this action is the $*$-algebra $k(\G,C(\bD))$
of functions $\phi : \G \to C(\bD)$ of finite support, with multiplication and involution given by
\[
\phi * \psi (\g_0) = \displaystyle \sum_{\g\in\G}\phi(\g)\alpha_{\g}(\psi(\g^{-1}\g_0))
\quad \text{and} \quad
\phi^*(\g)=\alpha_{\g}(\phi(\g^{-1})^*).
\]

The full crossed product algebra $C(\bD) \rtimes \G$ is the completion of the algebraic crossed product in an appropriate norm  \cite{ped}.  There is a natural embedding of $C(\bD)$ into $C(\bD) \rtimes \G$ which maps an element $f\in C(\bD)$ to the constant function taking the value $f$ on $\G$. 
The identity element $\id$ of $C(\bD) \rtimes \G$ is then identified with the constant function $\id(\w)=1,\, \w\in\bD$.
There is a natural unitary representation $\pi : \G \to C(\bD) \rtimes \G$, where $\pi(\g)$ is the function taking the value $\id$ at $\g$ and $0$ otherwise.
It is convenient to denote $\pi(\g)$ simply by $\g$. Thus a typical element of the dense $*$-algebra $k(\G,C(\bD))$ can be written as a finite sum $\sum_{\g} f_{\g} \g$, where $f_{\g}\in C(\bD)$, $\g\in\G$.
The definition of the multiplication implies the covariance relation
\[
\alpha_{\g}(f) = \g f \g^{-1} \quad \text{for} \quad f \in C(\bD), \g \in \Gamma.
\]

If $d\in \D^1$, let $\Om(d)$ denote the clopen subset of $\bD$ corresponding to the set of all semi-geodesics with initial edge $d$ and initial vertex $o(d)$. The indicator function $p_d$ of the set $\Om(d)$ is continuous and so lies in $C(\bD)\subset C(\bD) \rtimes \G$.

For each $d\in \D^1$, the element $p_d$ is a projection in $C(\bD) \rtimes \G$ and therefore defines an equivalence class $[p_d]$ in $K_0(C(\bD) \rtimes \G)$. (The standard reference for the K-theory of $C^*$-algebras is \cite{bl}.)  It is important to observe that edges $d_1, d_2$ lying in the same $\G$-orbit give rise to equivalent projections $p_{d_1}, p_{d_2}$, because of the covariance relations in the crossed product algebra (c.f. Lemma \ref{Kequivalence} below). Therefore the equivalence class $[p_d]$ depends only on the directed edge $e=\G d\in E$.
Write $[e]=[p_d]$. In this way, each edge $e\in E$ gives rise to an element $[e]\in K_0(C(\bD) \rtimes \G)$.

\begin{figure}[htbp]
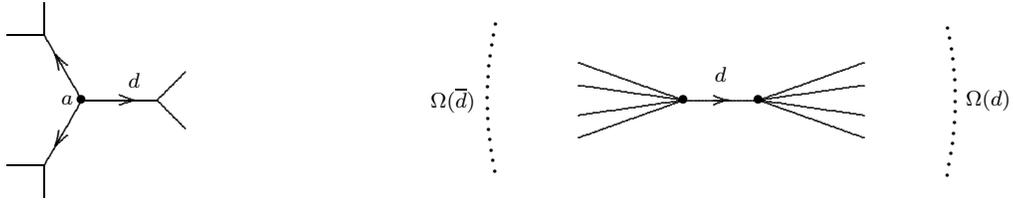
\label{relationdiag}
\hfil
\centerline{
\beginpicture
\setcoordinatesystem units <0.5cm,0.866cm>   
\setplotarea x from -4 to 4, y from -2  to 2         
\putrule from 0 0 to 2 0
\putrule from -2 1 to -1 1
\putrule from -2 -1 to -1 -1
\putrule from -1 1 to -1 1.5
\putrule from -1 -1 to -1 -1.5
\setlinear
\plot -1 1  0 0  -1 -1 /
\plot 2.75 0.433  2 0  2.75 -0.433  /
\put{$_\bullet$}  at  0 0
\put{$_d$}[b] at  1.4 0.2
\put{$_a$}[r] at  -0.2 0
\arrow <6pt> [.3,.67] from  1.2 0 to  1.4 0
\arrow <6pt> [.3,.67] from  -0.6 0.6 to  -0.7 0.7
\arrow <6pt> [.3,.67] from  -0.6 -0.6 to -0.7 -0.7
\setcoordinatesystem units <1cm, 1cm> point at -8 0  
\setplotarea  x from -4 to 4,  y from -2 to 3
\put{$_d$}   [b] at 0.5 0.25
\put{$_\bullet$}  at  0 0
\put{$_\bullet$}  at  1 0
\put{$_{\Omega(d)}$}   [l] at 3.75 0
\put{$_{\Omega(\od)}$}   [r] at -2.75 0
\arrow <6pt> [.2, .67] from  0.4 0  to 0.6 0
\setlinear
\plot   -1.4 0.5  0 0  1 0  2.4 0.5 /
\plot    1  0  2.4 0.2 /
\plot    1  0  2.4 -0.2 /
\plot    1  0  2.4 -0.5 /
\plot    0  0  -1.4 0.2 /
\plot    0  0  -1.4 -0.2 /
\plot    0  0  -1.4 -0.5 /
\setplotsymbol({$\cdot$}) \plotsymbolspacing=4pt
\circulararc 25 degrees from 3.5 -1 center at -1 0 
\circulararc 25 degrees from -2.5 1 center at 2 0
\endpicture
}
\hfil
\caption{The geometric meaning of relations (\ref{rel1}).}
\end{figure}

The projections $p_d$ satisfy the following relations as functions in $C(\bD)$, as illustrated in Figure~\ref{relationdiag}.
\begin{subequations}\label{rel1}
\begin{eqnarray}
\displaystyle\sum_{\substack{d\in \D^1 \\ o(d)=a}}p_d &=& \id  ,\qquad \text{for} \ a\in \D^0; \label{rel1a}\\
p_d + p_{\od} &=& \id   ,\qquad \text{for} \ d\in \D^1. \label{rel1b}
\end{eqnarray}
\end{subequations}
\refstepcounter{picture}
The relations (\ref{rel1}) project to the following relations in $K_0(C(\bD)~\rtimes~\Gamma)$.
\begin{subequations}\label{relK1}
\begin{eqnarray}
\displaystyle\sum_{\substack{e\in E \\ o(e)=v}}[e] &=& [\id]  ,\qquad \text{for} \ v\in V; \label{relK1a}\\
{[e]+[\ove]} &=& [\id]  ,\qquad \text{for} \ e \in E. \label{relK1b}
\end{eqnarray}
\end{subequations}
Since the map $e\mapsto o(e)\, :E\to V$ is surjective, the relations (\ref{relK1}) imply that 
\begin{equation*}
\begin{split}
n_0[\id]& = \sum_{v\in V}\displaystyle \sum_{\substack{e\in E \\ o(e)=v}}[e] = \sum_{e\in E} [e]\\
& = \sum_{e\in E_+}([e]+[\ove]) = \sum_{e\in E_+}[\id]\\
& = n_1[\id].
\end{split}
\end{equation*}
Therefore $(n_0-n_1).[\id]=0$.
This proves Theorem \ref{main}(1) in the $\tA_1$ case.

\medskip
\noindent The basic idea of this proof will be generalized to buildings of type $\tA_2$ in Section \ref{section3}.

\begin{remark}\label{treeremark}
If the tree $\D$ is homogeneous of degree $q+1$, with $q\ge 2$, then
$\chi(\G)=\frac{n_0}{2}(1-q)\ne 0$
and so the element $[\id]$ is torsion. 
\end{remark}
 
\bigskip

\section{The $\tA_n$ result}\label{section2}

Let $\D$ be a locally finite affine building of type~$\wt A_n$, where $n\ge 2$. 
Thus $\D$ is an $n$-dimensional simplicial complex whose $n$-simplices are called \emph{chambers}. Let $\D^k$ denote the set of $k$-simplices in $\D$ for $k=0,1,\dots,n$. 
Recall that each vertex $\dz\in\D^0$ is labelled with a type $\tau (\dz) \in \ZZ/(n+1)\ZZ$, and any chamber $\dii\in \D^n$ has precisely one vertex of each type. 
An \emph{apartment} in $\D$ is a subcomplex which is isomorphic to a Coxeter complex of type $\tA_n$. The building $\D$ is the union of its apartments. 

\refstepcounter{picture}
\begin{figure}[htbp]
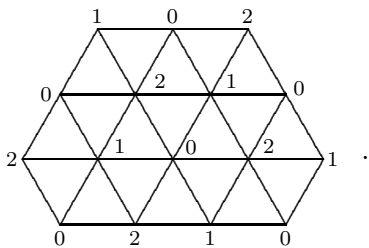
\label{A1}
\hfil
\centerline{
\beginpicture
\setcoordinatesystem units  <0.5cm, 0.866cm>        
\setplotarea x from -2.5 to 5, y from -2 to 2    
\put {$_1$} [b] at -2  2.1
\put {$_0$} [b] at  0  2.1
\put {$_2$} [b] at  2  2.1
\put {$_0$} [l] at 3.2 1.1
\put {$_1$} [l] at 4.1 0
\put {$_0$} [t] at 3  -1.1
\put {$_1$}  [t] at 1  -1.1
\put {$_2$} [t] at  -1  -1.1
\put {$_0$} [t] at -3 -1.1 
\put {$_2$} [r]  at  -4.1 0
\put {$_0$} [r] at -3.2  1
\put {$_2$} [l,b] at -0.5 1.1
\put {$_1$}  [l,b] at 1.4  1.1
\put {$_2$} [l] at  2.4 0.2
\put {$_1$}  [b]  at  -1.4 0.1
\put {$_0$} [t]  at  0.5  0.3
\putrule from -2   2     to  2  2
\putrule from  -3 1  to 3 1
\putrule from -4 0  to 4 0
\putrule from -3 -1  to  3 -1
\setlinear
\plot  -4 0  -2 2 /
\plot -3 -1  0 2 /
\plot -1 -1  2 2 /
\plot  1 -1  3 1 /
\plot  3 -1  4 0 /
\plot  -4 0  -3 -1 /
\plot  -3 1   -1 -1 /
\plot  -2 2   1 -1 /
\plot  0 2   3 -1 /
\plot  2 2   4 0 /
\endpicture.
}
\hfil
\caption{Part of an apartment in a building of type $\tA_2$, showing vertex types.}
\end{figure}

A \emph{sector} is a simplicial cone made up of chambers in some apartment and containing a unique base chamber \cite[Chapter 9]{ron}. Two sectors are  \emph{equivalent} if their intersection contains a sector. 
The boundary $\bD$ of $\D$ is by definition the set of equivalence classes of sectors in $\D$. For any $\omega \in \bD$ and $\dz\in \D^0$ there is a unique sector $[\dz,\w)$ in the class $\w$ having base vertex $\dz$ \cite[Theorem 9.6]{ron}.
The boundary is a totally disconnected compact Hausdorff space \cite[Section 2]{cms}, \cite[Section 4]{ca2}. 

Let $\G$ be a torsion free discrete group acting cocompactly on $\D$ by type rotating automorphisms. Then $X=\G\backslash\D$ is a cell complex with universal covering $\D$. Let $X^k$ denote the set of $k$-cells in $X$ for $k=0,1,\dots,n$.
There is a natural induced action of $\G$ on $\bD$, and we can form the universal crossed product $C^*$-algebra $C(\bD)\rtimes \G$.

\begin{lemma}\label{inversionfree}
The group $\G$ acts freely without inversion on $\D$. That is, no element $\g\in \G\setminus\{1\}$ stabilizes a nonempty simplex of $\D$.
\end{lemma}

\begin{proof}
Let $\delta$ be a nonempty simplex of $\D$ and define 
\[
G_\delta = \{g\in Aut(\D) \, ;\, g\delta=\delta\}.
\]
Then $G_\delta$ is a compact subgroup of the locally compact group $Aut(\D)$.
(See, for example, the proof of Theorem 4.1 of \cite {gs}.) Since $\G$ is a discrete subgroup of the automorphism group of $\D$, it follows that $\G\cap G_\delta$ is a finite subgroup
of $\G$. Since $\G$ is torsion free, $\G\cap G_\delta=\{ 1\}$.
\end{proof}

The first step is to formalize the notion of a directed cell in $X$. Let $\s^n$ be a model typed $n$-simplex with vertices $\0, \1, \2, \dots,\n$. Assume that the vertex $\p$ of $\s^n$ has type $\tau(\p)=p \in \ZZ / (n+1)\ZZ$. 

\refstepcounter{picture}
\begin{figure}[htbp]
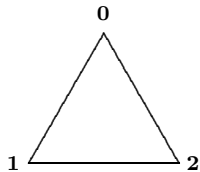
\label{modelsimplex}
\hfil
\centerline{
\beginpicture
\setcoordinatesystem units  <1cm, 1.732cm>
\setplotarea x from -1 to 1, y from -1 to 0.1         
\put {$_{\1}$} [r] at -1.1 -1
\put {$_{\2}$} [l] at 1.1 -1
\put {$_{\0}$} [b] at 0 0.1
\putrule from -1 -1 to 1 -1
\setlinear \plot -1 -1 0 0 1 -1 /
\endpicture
}
\hfil
\caption{The case $n=2$: the model simplex $\s^2$.}
\end{figure}

An isometry $r: \s^n \to \D$ is said to be {\em type rotating} if there exists $j\in \ZZ / (n+1)\ZZ$ such that, for each vertex $\p$ of $\s^n$, $\tau(r(\p))=\tau(\p)+j \pmod {n+1}$ . That is, $\tau(r(\p))=p + \tau(r(\0)) \pmod {n+1}$.
Let  $D^n$ denote the set of type rotating isometries $r: \s^n \to \D$.
Recall that any chamber $\dii\in \D^n$ has precisely one vertex of each type. Therefore, for each $j\in \ZZ / (n+1)\ZZ$ there is a unique $r \in D^n$ such that $r(\s^n)=\dii$ and $r(\0)$ has type $j$. An element $r\in D^n$ is therefore uniquely determined by the pair $(r(\s^n),r(\0))$. We therefore identify $D^n$ with the set of directed chambers of $\D$. Each geometric chamber $\dii \in \D^n$ is the image of each of the (n+1) directed chambers $\{ r\in D^n\ ; \ r(\s^n)=\dii\}$ under the map $r\mapsto r(\s^n)$.

\refstepcounter{picture}
\begin{figure}[htbp]
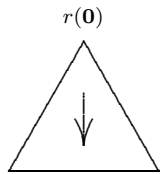
\label{directedchamber}
\hfil
\centerline{
\beginpicture
\setcoordinatesystem units  <1cm, 1.732cm>
\setplotarea x from -1 to 1, y from -1 to 0.1         
\arrow <10pt> [.2, .67] from  0 -0.4 to 0 -0.8
\put {$_{r(\0)}$} [b] at 0 0.1
\putrule from -1 -1 to 1 -1
\setlinear \plot -1 -1 0 0 1 -1 /
\endpicture
}
\hfil
\caption{The image of a directed chamber $r$.}
\end{figure}

The next lemma records this observation.

\begin{lemma}\label{nto1}
The map $r\mapsto r(\s^n)$ from $D^n$ to $\D^n$ is $(n+1)$-to-$1$.
\end{lemma}

Let $\fD^n = \G\backslash D^n$. Call $\fD^n$  the set of {\em directed $n$-cells} of $X=\G\backslash\D$.
There is a commutative diagram
\begin{equation*}
\begin{CD}
D^n                     @>r\mapsto r(\s^n)>>   \D^n\\
@VVV                        @VVV \\
\fD^n                  @>\e_n>>      X^n
\end{CD}
\end{equation*}
where the vertical arrows represent quotient maps and  $\e_n$ is defined by $\e_n(\G r)=\G.r(\s^n)$.

The next result makes precise the fact that each geometric $n$-cell in $X^n$ corresponds to exactly $n+1$ directed $n$-cells.

\begin{lemma}\label{epsilon}
The map $\e_n : \fD^n \to X^n$ is surjective and $(n+1)$-to-$1$.
\end{lemma}

\begin{proof}
Fix $\dii\in \D^n$.
By Lemma \ref{nto1}, $\{ r\in D^n\ ; \ r(\s^n)=\dii\}=\{r_0,r_1,\dots, r_n\}$, a set containing precisely $n+1$ elements.
Since $\G$ acts freely without inversion (Lemma \ref{inversionfree}), the set $\{\G r_0,\G r_1,\dots, \G r_n\} \subset \fD^n$ also contains precisely $n+1$ elements, each of which maps to $\G\dii$ under $\e_n$.

Now suppose that $\e_n(\G r)=\G\dii$ for some $r\in D^n$.
Then $\g r(\s^n)=\dii$ for some $\g\in\G$. Thus $\g r\in \{r_0,r_1,\dots, r_n\}$ and so
\[
\G r \in \{\G r_0,\G r_1,\dots, \G r_n\}.
\]
 This proves that $\e_n$ is $(n+1)$-to-$1$.
\end{proof}

Let $r\in D^n$; that is, suppose that $r: \s^n \to \D$ is a type
rotating isometry. Let
\begin{equation*}
\Omega(r)=\left \{ \omega \in \bD \ ;\ r(\s^n) \subset [r(\0),\omega) \right \},
\end{equation*}
the set of boundary points represented by sectors which originate at $r(\0)$ and contain $r(\s^n)$. (Figure \ref{Omegar} illustrates the case $n=2$.)

\refstepcounter{picture}
\begin{figure}[htbp]\label{Omegar}
\hfil
\centerline{
\beginpicture
\setcoordinatesystem units  <1cm, 1.732cm>
\setplotarea x from -1 to 1, y from -2 to 0.1         
\put {$_{r(\s^n)}$} [b] at 0 -0.8
\put {$_{r(\0)}$}  [b]  at 0 0.1
\put {$_{[r(\0),\w)}$} at 0 -2
\putrule from -1 -1 to 1 -1
\setlinear \plot -2 -2 0 0 2 -2 /
\endpicture
}
\hfil
\caption{The sector $[r(\0),\w)$, for $\w\in\Omega(r)$.}
\end{figure}

Clearly  $\Omega(\g r)=\g\Omega(r)$.
Let $p_r \in C(\bD)$ be the characteristic function of $\Om(r)$ and let $\id$ be the constant function $\id(\w)=1$, for $\w\in\bD$. 

The discussion in \cite[Section 2]{cms}, \cite[Section 4]{ca2} shows that for each vertex $\dz\in \D^0$,
\begin{equation}\label{1}
\id =  \displaystyle \sum_
{
\substack{
r\in D^n\\
r(\0)=\dz
}
}
p_r.
\end{equation}
The essential point is that each $\w\in\bD$ corresponds to a unique sector $[\dz,\w)$.
It follows from (\ref{1}) that in $K_0(C(\bD)\rtimes\G)$,
\begin{equation*}
[\id] =  \displaystyle \sum_
{
\substack{
r\in D^n\\
r(\0)=\dz
}
}
[p_r].
\end{equation*}
As in the $\tA_1$ case, it is important that the class $[p_r]$ in $K_0$ depends only on the orbit $\G r$.

\begin{lemma}\label{Kequivalence}
If $r, s \in D^n$ with $\G r=\G s$ then $[p_r]=[p_s]$.
\end{lemma}

\begin{proof}
If $r=\g s$ with $\g\in \G$ then the covariance condition for the action of $\G$ on $C(\Om)$ implies that $p_r=\g p_s\g^{-1}$. The result now follows because equivalent idempotents belong to the same class in $K_0$.
\end{proof}

In view of Lemma \ref{Kequivalence}, the following notation makes sense.

\begin{definition}
If $d=\G r \in \fD^n$, let $[d]=[p_r]\in K_0(C(\bD)\rtimes\G)$.
Also let $d(\p)=\G r(\p)\in X^0$ for $\p \in \{\0, \1,\dots,\n\}$, and let $d(\s^n)=\G r(\s^n)\in X^n$.
\end{definition}

Using this notation, for $x^0=\G\dz\in X^0$, 
\begin{equation*}
\{ [p_r]\ ;\ r\in D^n, r(\0)=\dz\}
= \{ [d]\ ;\ d\in \fD^n, d(\0)=x^0\}.
\end{equation*}

Therefore in $K_0(C(\bD)\rtimes\G)$,
\begin{equation}\label{2}
[\id] =  \displaystyle \sum_
{
\substack{
d\in \fD^n\\
d(\0)=x^0
}
}
[d], \qquad \text{for} \ \dz\in \D^0.
\end{equation}

Since the building $\D$ is of type $\tA_n$, where $n\ge 2$, it has additional structure. The link of any fixed vertex $\dz\in\D$ is the flag complex of a finite projective geometry $\Pi(\dz)$ of dimension $n$ and order $q\ge 2$ \cite[Theorem~3.5 and p.~95]{ron}.
Denote by $\Star(\dz)$ the set of vertices at distance~1 from $\dz\in\D$. There is a bijective correspondence between the vertices in $\Star(\dz)$ and subspaces of the projective geometry $\Pi(\dz)$. Two elements in $\Star(\dz)$ are incident in $\Pi(\dz)$ if they lie in a common chamber of $\D$.
If $n\ge 3$ then $\Pi(\dz)$ is a 
Desarguesian projective geometry $PG(n,q)$ \cite[(4.9)]{st} and $q$ is a prime power. On the other hand if $n=2$, then $\Pi(\dz)$ is a (possibly non-Desarguesian) projective plane. In either case, the order $q$ of $\Pi(\dz)$ is independent of $\dz$ and each $(n-1)$-dimensional simplex of $\D$ lies in $q+1$ chambers. The dimension $\Dim (x)$ of a subspace $x\in \Pi(\dz)$ is the maximum $k$ for which there is a flag of subspaces 
$x_1\subset x_2\subset\dots\subset x_k=x$. Thus a point of $\Pi(\dz)$ has dimension one and a line has dimension 2. This terminology has the unfortunate consequence that a maximal proper subspace of $n$-dimensional projective space has dimension $n$.
Note that if $x\in \Star (\dz)$ then $\tau (x)=\tau(\dz)+\Dim (x)$. 

Let $r\in D^n$ and fix $k\in \{1,2,\dots,n\}$.
Define, for each $r\in D^n$, 
\[
\fS_k(r)=\{ s\in D^n\ ;\ \Om(s)\subset\Om(r) \; \text{and} \;  s(\0)=r(\ka)\}.
\]

\refstepcounter{picture}
\begin{figure}[htbp]\label{transition}
\hfil
\centerline{
\beginpicture
\setcoordinatesystem units  <1cm, 1.732cm>
\setplotarea x from -1 to 1, y from -3 to 0.1         
\put {$_{r(\s^n)}$} [b] at 0 -0.8
\put {$_{s(\s^n)}$} [b] at -1 -1.8
\put {$_{r(\0)}$}  [b]  at 0 0.1
\put {$_{s(\0)=r(\ka)}$}  [r]  at -1.1 -1
\put {$_{s(\ka)}$}  [r]  at -2.1 -2
\put {$_{[r(\ka),\w)}$} at -1 -2.7
\putrule from -1 -1 to 1 -1
\putrule from -2 -2 to 0 -2
\setlinear \plot -2 -2 0 0 1 -1 /
\setlinear \plot -1 -1  0 -2 /
\setdashes \plot -2.5 -2.5  -2 -2 /
\setdashes \plot 0 -2   0.5 -2.5 /
\endpicture
}
\hfil
\caption{The condition $s\in\fS_k(r)$. }
\end{figure}
The set $\fS_k(r)$ contains $q^{k(n+1-k)}$ elements $s\in~D^n$.
In fact, a more general result is contained in \cite{ca2}. See in particular  Lemma 2.1 and the proof of Lemma 2.4 of that article.
We sketch the proof. 

Consider the set of chambers~$C$ such that
some sector based at $r(\mathbf{0})$ contains a subsector, based at
$r(\mathbf{k})$, whose initial chamber is~$C$. The set of such chambers~$C$ is
exactly the preimage of a
certain chamber with respect to the retraction from $r(\sigma^n)$.  As
such, the number of possible~$C$ is $q^l$ where $l$ is the number of
transitions in any stretched gallery between $r(\sigma^n)$ and any
such chamber~$C$.

As both of those chambers contain the vertex $r(\mathbf{k})$, one can
calculate the number of transitions working inside the spherical building
$\text{star}(r(\mathbf{k}))$.  In that spherical building, the chambers of
an apartment can be parameterized by the elements of the Weyl group,
which is just the symmetric group on $\{1,2,\dots,n+1\}$.  Fixing
specific coordinates one can identify $r(\sigma^n)$ with the identity permutation and $C$ with the permutation
\[
(k+1\quad k+2 \quad \dots \quad n+1 \quad 1 \quad 2 \quad \dots \quad k)
\]
The number of transitions in the stretched gallery we are considering
is equal to the length of this element in the Weyl group (relative to
the standard Coxeter basis, $(1 2)$, $(2 3)$, etc.), which in turn is
equal to the number of pairs whose order is inverted.  This last is
easily seen to be $k(n+1-k)$.

Thus $\fS_k(r)$ contains $q^{k(n+1-k)}$ elements $s\in~D^n$ and $\Om(r)$
 may be expressed as a disjoint union
\begin{equation*}
\Om(r)=  \displaystyle \bigsqcup_{s\in \fS_k(r)}\Om(s).
\end{equation*}
Therefore
\begin{equation}\label{sum0}
p_r =  \displaystyle \sum_{s\in \fS_k(r)} p_s.
\end{equation}
Moreover a dual argument shows that each $s\in D^n$ lies in precisely $q^{k(n+1-k)}$ of the sets $\fS_k(r)$, with $r\in D^n$.

Define as follows a \zomatrix\ $M_k$, with entries indexed by $\fD^n\times\fD^n$. If $c, d \in \fD^n$, then $M_k(d,c)=1$ if and only if there are representative isometries $r,s\in D^n$, with
$c=\G r$, $d= \G s$  and $s\in\fS_k(r)$; otherwise $M_k(d,c)=0$. Each row or column of $M_k$ has precisely $q^{k(n+1-k)}$ nonzero entries and  equation (\ref{sum0}) implies that,
for each $c\in \fD$, 
\begin{equation}\label{sum}
[c]= \displaystyle \sum_{d\in \fD^n} M_k(d,c) [d].
\end{equation}
Recall that $n_0$ denotes the number of vertices of $X=\G\backslash\D$, that is the number of $\G$-orbits of vertices of $\D$. It follows from equations (\ref{2}), (\ref{sum}) that
\begin{equation*}
\begin{split}
n_0[\id]& = \displaystyle \sum_{c\in \fD^n}[c] =
\displaystyle \sum_{c\in \fD^n}
\displaystyle \sum_{d\in \fD^n} M_k(d,c) [d] \\
& = \displaystyle \sum_{d\in \fD^n}\displaystyle \sum_{c\in \fD^n}M_k(d,c) [d]
= \displaystyle \sum_{d\in \fD^n}q^{k(n+1-k)} [d] \\
& = n_0q^{k(n+1-k)}[\id].
\end{split}
\end{equation*}
Therefore, for each $k\in \{1,2,\dots,n\}$,
\begin{equation}
n_0(q^{k(n+1-k)}-1)[\id] = 0.
\end{equation}
In particular, setting  $k=1,2$ gives
\[
n_0(q^n-1)[\id] = 0\qquad  \text{and} \qquad n_0(q^{2(n-1)}-1)[\id] = 0.
\]

\noindent If $n$ is odd then $\gcd(n,2(n-1))=1$ and so $n_0(q-1)[\id]= 0$.

\noindent If $n$ is even then $\gcd(n,2(n-1))=2$ and so $n_0(q^2-1)[\id]=0$.

This completes the proof of Theorem \ref{main}(2).\qed

\begin{remark}
The same argument, greatly simplified, applies in the case of a homogeneous tree $\D$ of degree $q+1$, where $q\ge 2$, and shows that 
$n_0(q-1)[\id]= 0$. This is weaker than the result obtained in Section \ref{section1}, Remark \ref{treeremark}.
\end{remark}

\bigskip


\section{The $\tA_2$ result}\label{section3}

This section will be concerned with buildings of type $\tA_2$. The previous results all apply, with $n=2$. The terminology and results from Section \ref{section2} will therefore be used without further comment. 

Let $\s^e$ denote the face $\{\1, \2\}$ of $\s^2=\{\0, \1, \2 \}$. We regard $\s^e$ as a model $1$-simplex. 
An isometry $r: \s^e \to \D$ is said to be {\em type rotating} if there exists $j\in \ZZ / 3\ZZ$ such that, for each vertex $\p$ of $\s^2$, $\tau(r(\p))=\tau(\p)+j \pmod 3$ for each vertex $\p$ of $\s^2$. That is, $\tau(r(\p))=p + \tau(r(\0)) \pmod 3$.
Let  $D^e$ denote the set of type rotating isometries $r: \s^e \to \D$.

\begin{figure}[htbp]
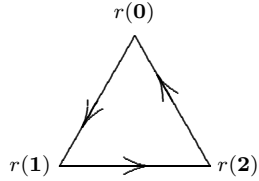
\label{edges}
\hfil
\centerline{
\beginpicture
\setcoordinatesystem units  <1cm, 1.732cm>
\setplotarea x from -1 to 1, y from -1 to 0.1         
\arrow <10pt> [.2, .67] from   0 -1 to 0.2 -1
\arrow <10pt> [.2, .67] from  -0.5 -0.5 to -0.7 -0.7 
\arrow <10pt> [.2, .67] from  0.5 -0.5  to 0.3 -0.3 
\put {$_{r(\1)}$} [r] at -1.1 -1
\put {$_{r(\2)}$} [l] at 1.1 -1
\put {$_{r(\0)}$} [b] at 0 0.1
\putrule from -1 -1 to 1 -1
\setlinear \plot -1 -1 0 0 1 -1 /
\endpicture
}
\hfil
\caption{Directed edges of a $2$-simplex in $\D$.}
\end{figure}

Any  $1$-simplex $\di \in \D^1$ has two vertices of two different types. The type rotating condition implies that there is a unique $r \in D^e$ such that $r(\s^e)=\di$. That is, $r\in D^e$ is uniquely determined by its image $r(\s^e)$.
We use $D^e$ to choose a preferred orientation for each edge of $\D$.
By a {\em directed edge}, we mean an edge with its preferred orientation. 
We therefore identify $D^e$ with the set of directed edges of $\D$.
Each geometric edge $\di \in \D^1$ is the image of a {\em unique} directed edge $r \in D^e$ under the map $r\mapsto r(\s^e)$. That is, each edge of $\D$ has a unique direction determined by the typing of the vertex set of $\D$. See Figure~\ref{edges}.  This is only superficially different from the $\tA_1$ case in Section \ref{section1}, where each geometric edge corresponds to two directed edges. In the tree case an edge is a maximal simplex.

\refstepcounter{picture}

\begin{lemma}
The map $r\mapsto r(\s^e)$ from $D^e$ to $\D^1$ is bijective.
\end{lemma}

Let $\fD^e = \G\backslash D^e$. Call $\fD^e$  the set of {\em directed edges} of $X=\G\backslash\D$.
There is a commutative diagram
\begin{equation*}
\begin{CD}
D^e                     @>r\mapsto r(\s^e)>>   \D^1\\
@VVV                        @VVV \\
\fD^e                  @>\e_1>>      X^1
\end{CD}
\end{equation*}
where the vertical arrows represent quotient maps and  $\e_1$ is defined by $\e_1(\G r)=\G.r(\s^e)$.

The next result makes precise the fact that each geometric edge in $X^1$ corresponds to exactly one directed edge (i.e. has a unique direction). It is an easy consequence of the previous lemma.

\begin{lemma}\label{node}
The map $\e_1 : \fD^e \to X^1$ is bijective.
\end{lemma}

If $r\in D^2$, let $\overline\Om(r)$ denote the set of $\w\in \bD$ such that $r(\0)$
is in the interior of the  convex hull $\conv(r(\s^n)\cup [r(\0),\w))$ as illustrated in Figure \ref{overlineOm}. It is easy to see that $\overline\Om(r)$ is a clopen set.

\refstepcounter{picture}
\begin{figure}[htbp]\label{overlineOm}
\hfil
\centerline{
\beginpicture
\setcoordinatesystem units <0.5cm,0.866cm>  
\setplotarea x from -4 to 4, y from 2 to 4   
\setlinear \plot  -1.5 3.5  1 1  3.5 3.5 /
\put {$_{[r(\0),\w)}$} [] at 0 3.7
\put{$_{\bullet}$} at 0 2
\put {$_{r(\0)}$} [] at 1  2
\putrule from -1 1 to 1 1
\setlinear \plot  -1 1    0 2  1.5  3.5 / 
\setlinear \plot -3.5 3.5  -1 1 /
\endpicture
}
\hfil
\caption{}
\end{figure}

Let $\op_r \in C(\bD)$ be the characteristic function of $\overline\Om(r)$. We have the following analogue of Lemma \ref{Kequivalence}, which has the same proof.

\begin{lemma}\label{Kequivalence'}
If $r, s \in D^2$ with $\G r=\G s$ then $[\op_r]=[\op_s]$.
\end{lemma}

In view of Lemma \ref{Kequivalence'}, the following notation makes sense.

\begin{definition}
If $d=\G r \in \fD^2$, let $[\od]=[\op_r]\in K_0(C(\bD)\rtimes\G)$.
\end{definition}

Recall now the following result.

\begin{lemma}\label{ronan}
{\rm \cite[Lemma 9.4]{ron}}
Given any chamber $\dd\in\D^2$ and any sector $S$ in $\D$, there exists a sector $S_1\subset S$ such that $S_1$ and $\dd$ lie in a common apartment.
\end{lemma}

\begin{figure}[htbp]
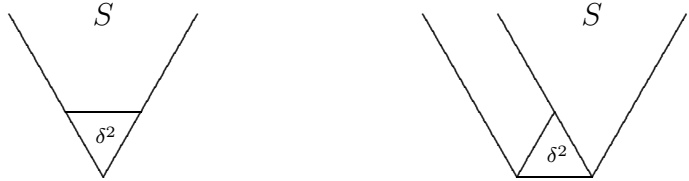
\label{chambersector}
\hfil
\centerline{
\beginpicture
\setcoordinatesystem units <0.5cm,0.866cm>    
\setplotarea x from -4 to 4, y from 1 to 3         
\put {$S$} [] at 0 2.5
\put {$_{\dd}$} [] at 0.1 0.65
\putrule from -1 1 to 1 1
\setlinear \plot -2.5 2.5  0 0  2.5 2.5 /
\setcoordinatesystem units <0.5cm,0.866cm>  point at -12 1 
\setplotarea x from -4 to 4, y from 0 to 4         
\setlinear \plot  -1.5 3.5  1 1  3.5 3.5 /
\put {$S$} [] at 1 3.5
\put {$_{\dd}$} [] at 0.1 1.35
\putrule from -1 1 to 1 1
\setlinear \plot  -1 1    0 2  / 
\setlinear \plot -3.5 3.5  -1 1 /
\endpicture
}
\hfil
\caption{Relative positions of a chamber and a representative sector.}
\end{figure}

It follows by considering parallel sectors in an appropriate apartment that 
if $\dd\in\D^2$ and if $\w\in\bD$, then $\w$ has a representative sector  $S$ that lies relative to $\dd$ in one of the two positions in Figure \ref{chambersector}, in some apartment containing them both.

\refstepcounter{picture}

Therefore 
\[
\id =  \displaystyle \sum_
{
\substack{
r\in D^2\\
r(\s^2)=\dd
}
}
p_r
+ \displaystyle \sum_
{
\substack{
r\in D^2\\
r(\s^2)=\dd
}
}
\op_r
\]
where there are three terms in each sum.

Now for $x^2=\G\dd\in X^2$, 
\[
\{ [p_r]\ ;\ r\in D^2, r(\s^2)=\dd\}
= \{ [d]\ ;\ d\in \fD^2, d(\s^2)=x^2\}.
\]
Therefore in $K_0(C(\bD)\rtimes\G)$,
\begin{equation}\label{3}
[\id] =  \displaystyle \sum_
{
\substack{
d\in \fD^2\\
d(\s^2)=x^2
}
}
[d]
+ \displaystyle \sum_
{
\substack{
d\in \fD^2\\
d(\s^2)=x^2
}
}
[\od].
\end{equation}

\bigskip

\begin{figure}[htbp]
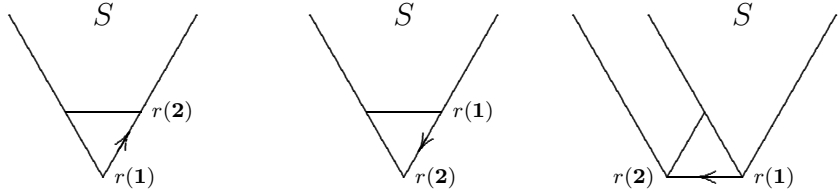
\label{edgesector}
\hfil
\centerline{
\beginpicture
\setcoordinatesystem units <0.5cm,0.866cm>    
\setplotarea x from -4 to 4, y from 1 to 3         
\put {$S$} [] at 0 2.5
\put {$_{r(\1)}$} [l] at 0.3 0
\put {$_{r(\2)}$} [l] at 1.3 1
\arrow <7pt> [.2, .67] from  0.4 0.4  to 0.7 0.7 
\putrule from -1 1 to 1 1
\setlinear \plot -2.5 2.5  0 0  2.5 2.5 /
\setcoordinatesystem units <0.5cm,0.866cm> point at -8 0   
\setplotarea x from -4 to 4, y from 1 to 3         
\put {$S$} [] at 0 2.5
\put {$_{r(\2)}$} [l] at 0.3 0
\put {$_{r(\1)}$} [l] at 1.3 1
\arrow <7pt> [.2, .67] from  0.6 0.6  to 0.4 0.4 
\putrule from -1 1 to 1 1
\setlinear \plot -2.5 2.5  0 0  2.5 2.5 /
\setcoordinatesystem units <0.5cm,0.866cm>  point at -16 1 
\setplotarea x from -4 to 4, y from 0 to 4         
\setlinear \plot  -1.5 3.5  1 1  3.5 3.5 /
\put {$S$} [] at 1 3.5
\put {$_{r(\2)}$} [r] at -1.3 1
\put {$_{r(\1)}$} [l] at 1.3 1
\arrow <7pt> [.2, .67] from  0 1  to -0.2 1 
\putrule from -1 1 to 1 1
\setlinear \plot  -1 1    0 2  / 
\setlinear \plot -3.5 3.5  -1 1 /
\endpicture
}
\hfil
\caption{Relative positions of an edge and a representative sector.}
\end{figure}

If $r\in D^e$ then any $\w \in \bD$ has a representative sector $S$ which lies relative to $r(\s^e)$ in one of the three positions in Figure \ref{edgesector}, in some apartment containing them both.
This is an easy consequence of Lemma \ref{ronan}.

\refstepcounter{picture}

Let $\Om(r), \overline\Om(r), \hat\Om(r)$ denote the corresponding clopen subsets of $\bD$, and let $q_r, \oq_r, \hq_r$ respectively denote the corresponding characteristic functions.  These functions are idempotents in $C(\bD)\rtimes\G$ and the equivalence classes 
$[q_r], [\oq_r], [\hq_r]$ in $K_0$ depend only on $\G r$. If $\di \in \D^1$ then there is a unique $r\in D^e$ such that $r(\s^e)=\di$ and we have $\id = q_r+\oq_r+\hq_r$.

\begin{definition}
If $e=\G r\in\fD^e$, where $r\in D^e$, let $[e]=[q_r]$, $[\ove]=[\oq_r]$ and
$[\he]=[\hq_r]$
\end{definition}

We now have
\begin{equation}\label{4}
[\id] = [e] + [\ove] + [\he]
\end{equation}

\begin{definition}
If $e=\G r\in\fD^e$ and $\p \in \{\1, \2\}$, let $e(\p)=\G r(\p)\in X^0$.
\end{definition}
\noindent Think of $e(\1)$ as the initial vertex and $e(\2)$ as the final vertex of the directed edge $e$. 

For each vertex $\dz\in \D^0$,
\begin{equation*}
\id =  \displaystyle \sum_
{
\substack{
r\in D^e\\
r(\1)=\dz
}
}
q_r
=  \displaystyle \sum_
{
\substack{
r\in D^e\\
r(\2)=\dz
}
}
\oq_r.
\end{equation*}
It follows that if $x^0=\G\dz\in X^0$ then
\begin{equation}\label{5}
[\id] =  \displaystyle \sum_
{
\substack{
e\in \fD^e\\
e(\1)=x^0
}
}
[e]
=  \displaystyle \sum_
{
\substack{
e\in \fD^e\\
e(\2)=x^0
}
}
[\ove].
\end{equation}
Given $r\in D^e$ and any chamber $\dd\in\D^2$ containing $r(\s^e)$, $r$ can
be extended uniquely to an element $s\in D^2$ such that $s(\s^2)=\dd$. Therefore
$\hat \Om(r)$ may be expressed as a disjoint union
\begin{equation*}
\hat \Om(r) = \bigsqcup \{\overline \Om(s)\; ; \; s\in D^2,\;  s\vert_{\s^e}=r \}.
\end{equation*}
It follows that $\hq_r=\sum\op_s$, where the sum is over $s\in D^2$ such that   $s\vert_{\s^e}=r$.
Choose one $r$ from each $\G$-orbit in $D^e$. Summing over these $r$ gives
\begin{equation*}
\displaystyle \sum_r \hq_r=
\displaystyle \sum_r
\displaystyle \sum_{s\mid_{\s^e}=r}
\op_s
\end{equation*}
where the double sum on the right contains precisely one $s$ from each $\G$-orbit in $D^2$,
because $\G$ acts freely without inversion on $\D$.
Taking equivalence classes in $K_0$ gives
\begin{equation}\label{6}
\displaystyle \sum_
{
e\in \fD^e
}
[\he]
=  \displaystyle \sum_
{
d\in \fD^2
}
[\od].
\end{equation}
\medskip

\noindent {\em Proof of Theorem \ref{main} in the $\tA_2$ case.}
\medskip

For $k=0,1,2$, let $n_k$ denote the number of $k$-cells in $X^k$.
By Lemma \ref{node},  $n_1=\#(\fD^e)$ and by Lemma \ref{epsilon}, $3n_2=\#(\fD^2)$.
It follows from equation (\ref{2}) with $n=2$, that
\begin{equation*}
n_0[\id] = \displaystyle \sum_{x^0\in X^0}
\displaystyle \sum_{\substack{d\in \fD^2\\d(\0)=x^0}}[d]
= \displaystyle \sum_{d\in \fD^2}[d].
\end{equation*}
Equation (\ref{3}) implies that
\begin{equation}\label{7'}
n_2[\id] = \displaystyle \sum_{x^2\in X^2}\left(
\displaystyle \sum_{\substack{d\in \fD^2\\d(\s^2)=x^2}}[d]
+ \displaystyle \sum_{\substack{d\in \fD^2\\d(\s^2)=x^2}}[\od]\right)
=\displaystyle \sum_{d\in \fD^2}[d]
+ \displaystyle \sum_{d\in \fD^2}[\od].
\end{equation}
Therefore
\begin{equation}\label{7}
(n_2-n_0)[\id]= \sum_{d\in \fD^2}[\od].
\end{equation}
Now by (\ref{5}),
\begin{equation}\label{9}
n_0[\id] =
\displaystyle \sum_{e\in \fD^e}[e]
=  \displaystyle \sum_{e\in \fD^e}[\ove].
\end{equation}
Using equations (\ref{4}),(\ref{9}),(\ref{6}) and (\ref{7}) successively yields
\begin{equation*}
\begin{split}
n_1[\id]& = \sum_{e\in \fD^e}[e] + \sum_{e\in \fD^e}[\ove] + \sum_{e\in \fD^e}[\he]\\
& = 2n_0[\id] + \sum_{e\in \fD^e}[\he]\\
& = 2n_0[\id] +  \displaystyle \sum_{d\in \fD^2}[\od]\\
& = 2n_0[\id] + (n_2-n_0)[\id]\\
& = n_0[\id] + n_2[\id].
\end{split}
\end{equation*}
This proves Theorem \ref{main}(1) in the $\tA_2$ case.
\qed

\bigskip

Note that the above proof has nowhere used the fact that the $\tA_2$ building $\D$ has order $q$.

\bigskip

\section{Euler-Poincar\'e characteristic and the $\tA_2$ case} \label{EPC}

In this section we sharpen the result of Theorem \ref{main} in the $\tA_2$ case. We assume throughout that $\D$ is a building of type $\tA_2$ and order $q$, and that $\G$ is a torsion free discrete subgroup of $\aut(\D)$ which acts cocompactly on $\D$.

The first step is to obtain an explicit expression for $\chi(\G)$.
Recall that $\chi(\G)$ coincides with the usual Euler-Poincar\'e characteristic of cell complex $X=\G\backslash\D$ \cite[I.1.5, Proposition 9]{ser2}.
That is $\chi(\G)=n_0-n_1+n_2$, where $n_k$ denotes the number of $k$-cells in $X^k$.

Since the building $\D$ is of type $\tA_2$, the link of any fixed vertex $\dz\in\D$ is the flag complex of a finite projective plane $\Pi(\dz)$ of order $q$. The vertices in $\Star(\dz)$ are the points and lines of $\Pi(\dz)$. A point and a line in $\Star(\dz)$ are incident in $\Pi(\dz)$ if they lie in a common chamber of $\D$.
Now $\Pi(\dz)$ contains $q^2+q+1$ points and $q^2+q+1$ lines. Moreover, each point is incident with $q+1$ lines. Therefore the vertex $\dz\in \D$ lies on $2(q^2+q+1)$ edges in $\D^1$ and on $(q+1)(q^2+q+1)$ chambers in $\D^2$. Since there are $n_0$ vertex orbits in $\D^0$, $\#(\fD^1)=2n_0(q^2+q+1)$ and $\#(\fD^2)=n_0(q+1)(q^2+q+1)$.
It follows from Lemma \ref{epsilon} that $n_1=n_0(q^2+q+1)$ and $n_2=n_0(q+1)(q^2+q+1)/3$.  Therefore

\begin{equation}\label{epmeas}
\chi(\G)=\frac{n_0}{3}(q-1)(q^2-1).
\end{equation}

\begin{remark}
Let $\FF$ be a non-archimedean local field with residue field of order $q$. Let $\G$ be a torsion free lattice in $G = \PGL(3,\FF)$. Recall that such a lattice is necessarily cocompact. 
The group $G$ acts in a type rotating manner on its Bruhat-Tits building $\D$ \cite{st}. The building $\D$ is of type $\tA_2$ and Theorem~\ref{main} therefore applies to any torsion free lattice $\G$ in $G$.  By \cite[3.3]{ser2}, $G$ has a unique invariant measure $\mu_G$ (called the {\em canonical measure} of $G$) which is an {\em Euler-Poincar\'e measure}. That is $\mu_G(G/\G)=\chi(\G)$, for each torsion free lattice $\G$ in $G$.

In the present situation it is not hard to determine $\mu_G$ directly. Let $\cO$ denote the valuation ring of $\FF$. Then $K=\PGL(3,\cO)$ is an open maximal compact subgroup of $G$ and the vertex set of $\D$ may be identified with the homogeneous space $G/K$. Then arguing as in \cite[p. 150, Th\'eor\`eme 7]{ser2} shows that
$\mu_G$ is the Haar measure on $G$ normalized so that
$\mu_G(K)=\frac{1}{3}(q-1)(q^2-1)$, in agreement with (\ref{epmeas}).
\end{remark}

The statement of Theorem \ref{main} can now be improved in the $\tA_2$ case.

\begin{corollary}\label{gcd}
Let $\G$ be a torsion free discrete group  of type rotating automorphisms of
a locally finite affine building $\D$ of type $\tA_2$ and order $q$. 
Suppose that $\G$ acts cocompactly on $\D$ with a finite number $n_0$ of vertex orbits. 
Then  $m.[\id]=0$ in the group $K_0(C(\bD)\rtimes \G)$, where
\begin{equation*}
m = 
\begin{cases}
    \frac{n_0}{3}(q^2-1)&  \text{if $q \not\equiv 1 \pmod {3}$},\\
    n_0(q^2-1) &  \text{if $q  \equiv 1 \pmod {3}$}.
\end{cases}
\end{equation*}
Moreover, if $q \equiv 0 \pmod {3}$, then $n_0$ is a multiple of $3$.
\end{corollary}

\begin{proof}
Theorem \ref{main} implies that $m.[\id] = 0$, where $m$ is the greatest common divisor of $n_0(q^2-1)$ and $\chi(\G)$. There are different outcomes, depending on the value of  $q \pmod {3}$.
\begin{enumerate}
\item  If $q \equiv 0 \pmod {3}$, then $n_0$ is a multiple of $3$, since 
$\chi(\G)$ is an integer. Using the fact that $(3,q-1)=1$ gives
\[
m=\big( 3\frac{n_0}{3}(q^2-1), (q-1)\frac{n_0}{3}(q^2-1) \big)=
\frac{n_0}{3}(q^2-1).
\]
\item If $q \equiv 1 \pmod {3}$ then  $\chi(\G)= \frac{q-1}{3}n_0(q^2-1)$ and so $m=n_0(q^2-1)$. This does not improve Theorem \ref{main}(2).
\item If $q \equiv 2 \pmod {3}$ then again the fact that $(3,q-1)=1$ implies that
\[
m=\big( 3n_0\frac{(q^2-1)}{3}, (q-1)n_0\frac{(q^2-1)}{3} \big)=
\frac{n_0}{3}(q^2-1).
\]
\end{enumerate}
\end{proof}

\begin{example}\label{a2exs}
Suppose  that $\G$ is a discrete subgroup of $\aut(\D)$ which acts freely and {\em transitively} in a type rotating manner on the vertex set of $\D$. That is, $n_0=1$. We refer to such a group $\G$ as an $\tA_2$ group \cite{cmsz}.
An $\tA_2$ group may have $3$-torsion and stabilize a chamber of $\D$. In fact this always happens if $q \equiv 0 \pmod {3}$.

Let $\G$ is a {\em torsion free} $\tA_2$ group.
If $q=2$ then there are five possible such groups $\G$, denoted A1, A2, B1, B2, C1 in \cite{cmsz}. These groups  are lattices in $\PGL(3,\FF)$, where $\FF= \FF_2((X))$
(the groups A1, A2) or $\FF=\QQ_2$ (the groups B1, B2, C1).
Corollary \ref{a2exs} implies that $[\id]=0$ in these cases, a fact confirmed by the computations in \cite{rs3}. 

Other examples of torsion free $\tA_2$ groups are the {\em regular} groups for $q=2,4,5,7,8,11$, constructed in \cite[I, Section 4]{cmsz}. These groups are lattices in $\PGL(3,\FF)$, where $\FF=\FF_q((X))$.
In particular, if $\G$ is the regular $\widetilde A_2$ group
with $q=5$, then $\frac{(q^2-1)}{3}= 8$ whereas direct calculation \cite{rs3} gives $4.[\id]=0$.
Therefore the estimates for the order of $[\id]$ given by Corollary \ref{gcd} are not best possible.
\end{example}

\begin{remark}
In \cite{rs3,rs4}, T. Steger and the author performed extensive computations which determined the order of $[\id]$ for  many $\tA_2$ groups with $q\le 11$. The computations were done for all the $\tA_2$ groups in the cases $q=2,3$ and for several representative groups for each of the other values of $q\le 11$.
If $q=2$ there are precisely eight $\tA_2$ groups $\G$, all of which embed as lattices in a linear group $\PGL (3,\FF)$ where $\FF= \FF_2((X))$ or $\FF=\QQ_2$. If $q=3$ there are 89 possible $\tA_2$ groups, of which 65 do not embed naturally in linear groups.

The results of the computations are indicated in the table below.

\bigskip

\centerline{
{\small
\begin{tabular}{|l|l|l|l|l|l|l|l|l|}
\hline
$q$                  & 2 & 3 & 4 & 5 & 7 & 8 & 9 & 11   \\ \hline
order of $[\id] \;$ & 1 & 2 & 1 & 4 & 2 & 7 & 8 & 10   \\ \hline
\end{tabular}
}
}

\bigskip

This experimental evidence suggests that for boundary crossed product algebras associated with $\tA_2$ groups it is always true that
$[\id]$ has order $q-1$ for $q \not\equiv 1 \pmod {3}$ and has order $(q-1)/3$ for $q  \equiv 1 \pmod {3}$.  It was proved in \cite[Proposition 8.3]{rs3} that the order of $[\id]$ cannot be smaller than these conjectured values.

It is striking that the order of $[\id]$ appears to depend {\em only} on the parameter $q$. Moreover the order is the same whether or not the group is torsion free. This suggests that the restriction to torsion free groups in the present article may eventually prove unnecessary.

It is tempting to conjecture that under the general hypotheses of Theorem \ref{main}, where $\G$ need not act transitively on the vertices of $\D$,  the order of $[\id]$ depends only on $q$ and the covolume $n_0$.
\end{remark}

\begin{remark}\label{EPCn}
If the building $\D$ is of type $\tA_n$, where $n\ge 3$, then according to a well known Theorem of J. Tits \cite[p. 137]{ron}, 
$\D$ is the building of $\PGL(n+1,\FF)$, where $\FF$ is a (possibly non-commutative) non-archimedean local field.

Let $\G$  be a torsion free lattice in $\PGL(n+1,\FF)$, where $n\ge 1$ and $\FF$ is a non-archimedean local field with residue field of order $q$. Then arguing as in \cite[p. 150, Th\'eor\`eme 7]{ser2} gives
\begin{equation}\label{epmeasn}
\chi(\G)=\frac{(-1)^n}{n+1}n_0(q-1)(q^2-1)\dots(q^n-1),
\end{equation}
where, as usual, $n_0$ is the number of vertices of $X=\G\backslash\D$.
If $n\ge 3$, the second part of Theorem \ref{main} does not imply that $\chi(\G).[\id]=0$. For example, if $q=n=4$ then $\chi(\G)$ is not divisible by $(q^2-1)$. It seems likely however that in general, $\chi(\G).[\id]=0$.
Note that it follows from \cite[XIII Theorem 2.6]{bw} that the Betti numbers of $\G$ satisfy 
$b_i(\G)=0$ for $i\ne 0,n$ and therefore $\chi(\G)=1+(-1)^nb_n(\G)$.
\end{remark}

\end{document}